\documentclass[12pt]{article}

\usepackage{bussproofs}
\EnableBpAbbreviations

\usepackage{url}

\def\HDS{\vrule width0pt height2.3ex depth1.05ex\displaystyle}
\def\f#1#2{{{\HDS #1}\over{\HDS #2}}}

\begin{document}

\title{{\LARGE G\" odel's Natural Deduction}}
\author{{\sc Kosta Do\v sen} and {\sc Milo\v s Ad\v zi\' c}
\\[1ex]
{\small Faculty of Philosophy, University of Belgrade}\\[-.5ex]
{\small \v Cika Ljubina 18-20, 11000 Belgrade, Serbia, and}\\
{\small Mathematical Institute, Serbian Academy of Sciences and Arts}\\[-.5ex]
{\small Knez Mihailova 36, p.f.\ 367, 11001 Belgrade, Serbia}\\[.5ex]
{\small email: kosta@mi.sanu.ac.rs, milos.adzic@gmail.com}}
\date{\small March 2016}
\maketitle

\begin{abstract}
\noindent This is a companion to a paper by the authors entitled ``G\" odel on deduction'', which examined the links between some philosophical views ascribed to G\" odel and general proof theory. When writing that other paper, the authors were not acquainted with a system of natural deduction that G\" odel presented with the help of Gentzen's sequents, which amounts to Ja\' skowski's natural deduction system of 1934, and which may be found in G\" odel's unpublished notes for the elementary logic course he gave in 1939 at the University of Notre Dame. Here one finds a presentation
of this system of G\" odel accompanied by a brief reexamination in the light of the notes of some points concerning his interest in sequents made in the preceding paper. This is preceded by a brief summary of G\" odel's Notre Dame course, and is followed by comments concerning G\" odel's natural deduction system.
\end{abstract}

\vspace{2ex}


\vspace{2ex}


\vspace{2ex}

\noindent {\small \emph{Acknowledgements.} Work on this paper was
supported by the Ministry of Education, Science and Technological Development of Serbia,
while the Alexander von Humboldt Foundation has supported the
presentation of matters related to it by the first-mentioned of us at the conference {\it General Proof Theory:
Celebrating 50 Years of Dag Prawitz's ``Natural Deduction''}, in
T\" ubingen, in November 2015. We are grateful to the organizers for the invitation to the
conference, for their care, and in particular to Peter Schroeder-Heister, for his exquisite hospitality.
We are also very grateful to Gabriella Crocco for making some of G\" odel's
unpublished notes available to us and for allowing a quotation out of them, which was made
readable through the decipherment she supervised.
The publishing of these notes is part of the project {\it Kurt G\" odel Philosopher: From
Logic to Cosmology}, which is directed by her and funded by the French
National Research Agency (project ANR-09-BLAN-0313). She and Antonio Piccolomini d'Aragona were also very kind to invite us to the workshop {\it Inferences and Proofs}, in Marseille, in May 2016, where the first-mentioned of us delivered a talk based partly on this paper, and enjoyed their especial hospitality.
We are grateful to the Institute for Advanced Study in Princeton for granting us, through the office of its
librarian Mrs Marcia Tucker and its archivist Mr Casey Westerman, the permission to
use G\" odel's unpublished notes; we were asked to give credit for that with the following text:
``All works of Kurt G\" odel used with permission. Unpublished Copyright
(1934-1978) Institute for Advanced Study. All rights reserved by Institute
for Advanced Study.''}

\vspace{5ex}

\section{Introduction}
Among G\" odel's unpublished writings at the Princeton University Library, two sets of notes for courses on elementary logic have been preserved (see \cite{Daw05}, Section 1.II, pp.\ 152-154). He gave the first of these courses at the University of Vienna in the summer of 1935 (see \cite{Daw97}, p.\ 108). According to \cite{Daw05} (p.\ 153), the notes for the Vienna course are about: ``\ldots truth tables, predicate logic, Skolem normal form, the Skolem-L\" owenheim theorem, the decision problem, and set theory''.

The second of these courses of G\" odel's was for graduate students at the University of Notre Dame in the spring semester of 1939 (see \cite{Daw97}, pp.\ 135-136, and \cite{Daw}). Among the topics G\" odel covered in these notes, which we will summarize in Section~3 below, we are here interested in particular in a natural deduction system for propositional logic, which G\" odel presented by relying on Gentzen's sequents. G\" odel does not use the terms \emph{natural deduction} and \emph{sequent}, but he does refer to Gentzen.

We considered in Section~2 of \cite{DA16} the interest these matters have for the understanding of the links between some philosophical views ascribed to G\" odel and general proof theory, which is the subject matter of that paper. When we wrote that, we were not acquainted with G\" odel's notes for the Notre Dame course, except for very brief summaries and a few fragments in Dawson's writings mentioned above and the bigger extracts in \cite{CN09}. In the meantime (see the Acknowledgements), we were able to see the relevant portion of G\" odel's notes, and our main goal is to present G\" odel's natural deduction system and make comments concerning it. Before doing that, in Section~3 we reexamine in the light of G\" odel's unpublished notes a few matters based on what we said in \cite{DA16} concerning G\" odel's interest in sequents. This will put our presentation of G\" odel's system within the context of general proof theory provided by that paper.

To put our presentation of G\" odel's natural deduction system in Section~4 within the context of his whole Notre Dame course, we summarize and make a few comments about his notes for this course in the next section. In Sections 5 and 6, we make comments upon this system.

\section{The Notre Dame course}
G\" odel gave at Notre Dame a one-semester course in elementary logic. His aim was to introduce his students to propositional and predicate logic, and also briefly at the end to the theory of types. We summarize the course in more detail, and make comments upon it, in \cite{AD16}. Short summaries may also be found in \cite{Daw97} (p.\ 136) and \cite{CN09} (Section~2), where there are also three larger extracts and comments upon them.

The course was influenced much by \cite{HA28}. As that book, it introduces the student to most of the major problems of modern logic, characteristic for its main branches. It stresses in particular the importance of the problem of completeness, and how the successful treatment of it makes a tremendous advance in comparison with what we had in the old Aristotelian tradition. It deals also with the problem of decidability, which is as much beyond the reach of the old logic, and the problem of the independence of axioms.

G\" odel explains with care truth-functional material implication, and tells that it corresponds as closely to \emph{if then} as a precise notion can correspond to an imprecise notion of ordinary language (see \cite{CN09}, p.\ 83, fn 14). Besides this truth-functional, extensional, implication, he envisages an intensional implication tied to deduction. G\" odel claims that the intensional implication in question here is the strict implication of modal logic, but intuitionistic implication might be better connected with deduction, in the light of \cite{G35}, and later developments in general proof theory.

G\" odel finds that the simple solution of the decidability problem for propositional logic is another great advance made by modern logic over traditional Aristotelian logic. This advance is made by espousing the truth-functional point of view. With an intensional implication, decidability is more difficult to obtain.

The formal system for propositional logic of the Hilbert type that G\" odel has in his notes is taken from \cite{HA28} (Section I.10). It is the system of \cite{WR10} (Section I.A.$\ast$1) without one axiom, which was found by Bernays to be provable from the remainder. The remaining axioms were shown independent again by Bernays, and G\" odel proves independence for one of the axioms, as an example. (At the time when he gave his course in Notre Dame, G\" odel worked on the independence problem in set theory for the Axiom of Choice and the Continuum Hypothesis.)
For this propositional axiom system, G\" odel proves soundness and later completeness, by proceeding as in \cite{Ka35}. G\" odel also proves functional completeness for the connectives of disjunction and negation, and shows that equivalence and negation are not functionally complete.

Since we deal with deduction in this paper, it is worth mentioning that G\" odel echoes the morale of \cite{Ca95} by noting that inference rules are indispensable for setting up a formal system, and cannot be replaced by implicational axioms.

Before turning to predicate logic, G\" odel makes remarks that we will survey in Section~3, and he deals with sequents and his natural deduction system, which will occupy us in Sections 4-6.

G\" odel's treatment of predicate logic is more cursory than his treatment of propositional logic. The completeness of first-order predicate logic is mentioned, but not proved. The same holds for the undecidability of predicate logic, and the decidability of monadic predicate logic.

G\" odel adds to the system of the Hilbert type he had for propositional logic one axiom and one rule involving the universal quantifier, both taken from \cite{HA28} (Section III.5), and adapts the substitution rule to cover individual variables, to get an axiom system for first-order predicate logic.

He uses the term ``tautology'' not only for the universally true formulae of propositional logic, but also for such formulae of predicate logic, which nowadays are called rather \emph{valid formulae}. He says however that using ``tautology'' is better abandoned if it is tied to the philosophical position that logic is devoid of content, that it says nothing. As the rest of mathematics, logic should be indifferent towards this position. (A citation covering this matter is in \cite{CN09}, p.\ 73.)

In making this comment, G\" odel might have had in mind Wittgenstein (who is usually credited for introducing the term ``tautology'' in modern logic, though the term was already used by Kant in a related sense).
The philosophical opinions of the early Wittgenstein, and perhaps also of the later, as well as those of the logical positivists, could be described by saying that they thought that logic is devoid of content.

G\" odel speaks of matters of completeness and decidability with a simile (see \cite{CN09}, Extract~2, pp.\ 84-85). Turing machines are not named, but their working is suggested by ``thinking machines''. One such device has a crank, which has to be turned to produce tautologies, where this word covers also the universally true formulae of predicate logic, and there is another device with a typewriter and a bell, which rings if one types in a tautology. The device of the second kind is available for propositional logic and monadic predicate logic, but not for the whole predicate calculus.

G\" odel rejects the unjustified existential presuppositions of Aristotelian syllogistic. This he does because they are either an empirical matter foreign to logic, or they would hamper arguments where they are not made, in which the issue might be exactly whether a predicate applies to anything or not.

The last part of the course deals with matters that lead to the theory of types. He considers classes as extensions of unary predicates, and briefly mentions Russell's ``no class theory'', which is referred to in \cite{DA16} (Section~4). Russell's convention for understanding definite descriptions is presented in a few sentences printed in \cite{CN09} (pp.\ 70-71), and Russell's achievement is succinctly presented as making the meaningfulness of language independent of empirical matters.

At the end of the course, G\" odel turns to paradoxes, of which he mentions the Burali-Forti paradox and Russell's paradox. The latter is treated in more detail, and there is a brief discussion of the theory of types. G\" odel's remarks on paradoxes are published as Extract~3 in \cite{CN09} (pp.\ 85-89).

\section{G\" odel's interest in sequents}
It is not well known that G\" odel had a favourable opinion about Gentzen's presentation of logic with sequents introduced in \cite{G35}. In the elementary logic course he gave at Notre Dame his aim was ``(1) to give, as far as possible, a complete theory of logical inferences and [of] logically true propositions, and (2) to show how they can be reduced to a certain number of primitive laws [all of them can be deduced from a minimum number of primitive laws]'' (this quotation from G\" odel's notes was printed in Extract~1 of \cite{CN09}, p.\ 78; we give here in square brackets also the version from Notebook 0). After proving completeness for the axiomatization of propositional logic of the Hilbert type, which we considered in the preceding section, G\" odel made some comments upon this result.

He wrote: ``I wish to stress that the interest of this result does not lie so much in this that our particular four axioms and three rules are sufficient to deduce everything, but the real interest consists in this that for the first time in the history of logic it has really been \emph{proved} that one \emph{can} reduce all laws of a certain part of logic to a few logical axioms. You know it has often been claimed that this can be done, and sometimes the law of contradiction and excluded middle have been considered as the logical axioms. But not even the shadow of a proof was given that really every logical inference can be deduced from them. Moreover, the assertion to be proved was not even clearly formulated, because it means nothing to say that something can be derived, e.g.\ from the law of contradiction, unless you specify in addition the rules of inference which are to be used in the derivation. As I said before, it is not so very important that just our four axioms are sufficient. After the method has once been developed, it is possible to give many other sorts of axioms which are also sufficient to derive all tautologies of the calculus of propositions. I have chosen the above four axioms because they are used in the standard textbooks of logistics. But I do not at all want to say that this choice was particularly fortunate. On the contrary, our system of axioms is open to some objections from the aesthetic point of view; e.g.\ one of the aesthetic requirements for a set of axioms is that the axioms should be as simple and evident as possible---in any case simpler than the theorems to be proved, whereas in our system e.g.\ the last axiom is pretty complicated, and on the other hand the very simple law of identity $p\supset p$ appears as a theorem. So in our system it happens sometimes that simpler propositions are proved more complicated than the axioms, which is to be avoided if possible. Recently, a system that avoids these disadvantages was set up by the mathematician G.\ Gentzen.''

The rendering of G\" odel's text here is not literal and scholarly. We have made slight, obvious, corrections to what he wrote, and added some punctuation marks. Parts of this text may be found in \cite{CN09} (p. 70) and \cite{Daw97} (p.\ 136). G\" odel's way of presenting logic with sequents is mentioned in the summary and comments preceding the three extracts from G\"odel's notes for the Notre Dame course published in \cite{CN09} (p.\ 70, to which we referred above), but the part of the notes with this presentation was not chosen there for printing.

We commented in Section~2 of \cite{DA16} upon the importance that G\" odel might have attached to the \emph{aesthetic} matters he mentions. We suggested that this was not a slight matter to him, and accords well with the platonism of his world view. It is not however something peculiar to G\" odel, but is also important for most, if not all, true mathematicians. We considered this matter more extensively towards the end of Section~2 of \cite{DA16}, and quoted also words ascribed to G\" odel that he should have said in his old age. Here we repeat only a short succinct statement, which we believe summarizes nicely G\" odel's opinion. It is taken from his unpublished notes written a few years after the Notre Dame course:

\begin{quote}
The truth is what has the simplest and the most beautiful \nopagebreak{symbolic expression.}
\end{quote}
(Acknowledgement for the translation of this sentence of \cite{GoX}, p.\ \textbf{[18]}, may be found in \cite{DA16}, Section~2.)

We considered in the remainder of Section~2 of \cite{DA16} whether it may be taken that G\" odel understood sequents as saying something about deductions, and not only as a peculiar way to write implications. G\" odel's natural deduction system may perhaps suggest that. The connection of sequents with natural deduction, from which they stem, and what Gentzen does with them, suggest strongly that they should be understood as being about deductions, as indeed many of those who worked with them in the last decades understand them. G\" odel's notes about his natural deduction system, which we will present in Section~4 below, do not give however a clear answer concerning that matter. There he calls sequents \emph{implications}, and often speaks as if the arrow that he writes in the middle of a sequent (instead of the turnstile, which is also sometimes written) was not essentially different from the horseshoe he writes for implication. We should note however that G\" odel was speaking to an audience of novices, and moreover the time had not yet come for the ideas of general proof theory that Gentzen suggested strongly, but that he himself did not espouse quite explicitly.

\section{The natural deduction system}
We believe that in 1939 it was quite unusual to present propositional logic in a logic course for beginners with a natural deduction system. The time for that had not yet come. It is in the second half of the last century that natural deduction entered this ground. When after the Second World War elementary textbooks with some kind of natural deduction started appearing in America (see \cite{PH12}), the dominant style was closer to \cite{J34} than to \cite{G35}, which was followed by \cite{P65}. G\" odel's own natural deduction system in the notes for the Notre Dame course can be briefly described as being made of Ja\' skowski's rules presented with Gentzen's sequents---Ja\' skowski dressed in Gentzen's garb. (G\" odel mentions Gentzen, but not Ja\' skowski.)

The presentation of natural deduction with sequents is implicit already in \cite{G35}, and is quite explicit in \cite{G36} (Section II.5). It is probable that G\" odel was inspired by that later paper, which he knew and cited in 1938 (see the beginning of Section~2 of \cite{DA16}). Gentzen's rules are however, as usual in natural deduction, \emph{multiplicative}, to use the modern terminology of substructural logics, while G\" odel's implication elimination is \emph{additive}. Gentzen assumes the structural rules of permutation, contraction and thinning, while G\" odel has only two sorts of thinning (see the comments concerning structural rules in Section~5 below). Their rules for negation are not the same.

G\" odel's propositional language in his notes for the course has, in the notation he uses, which we follow throughout, the propositional letters $p,q,r,\ldots$ and the connectives of negation $\sim$ and implication $\supset$. He uses the capital letters $P,Q,R,\ldots$ as schematic letters for the propositional formulae of the language with these primitive connectives, which he will call later \emph{expressions of the first kind}. He uses the capital Greek letters $\Delta,\Gamma,\ldots$ as schematic letters for sequences of these propositional formulae, presumably finite (though he says there are ``an arbitrary number'' of them) and possibly empty. As Gentzen, G\" odel conceives of sequents, which he calls ``secondary formulas'', as being also words of a language. This is why in a sequent $\Delta\rightarrow P$, the sequence of propositional formulae $\Delta$ should be finite. G\" odel calls here $\rightarrow$, which is often written $\vdash$, but which Gentzen too writes $\rightarrow$, ``another kind of implication''. Since he does not use the term \emph{sequent}, as we said in the Introduction (Section~1), this corroborates up to a point our supposition from \cite{DA16} (Section~3, concerning 8.4.16, and Section~5, concerning 8.4.18 [part I]) that in the 1970s he could have called sequents \emph{implications}. G\" odel also says in the Notre Dame notes that $\Delta\rightarrow P$ is a formula, but of another kind than $P$, and interprets it as saying that the sequence $\Delta$ of expressions of the first kind, i.e.\ sequence of propositional formulae, implies the expression of the first kind, i.e.\ propositional formula, $P$.

Then he presents the axioms and rules of inference for his natural deduction system. As axioms he has only the identity law, i.e.\ all sequents of the form $P\rightarrow P$, and as rules of inference he has the structural rule of thinning, which he calls ``addition of premisses'':
\[
\f{\Delta\rightarrow Q}{P,\Delta\rightarrow Q}\quad\quad\quad\f{\Delta\rightarrow Q}{\Delta, P\rightarrow Q}
\]
the rules for introducing and eliminating implication, which he calls ``the rule of exportation'' and ``the rule of implication'' respectively:
\[
\f{\Delta, P\rightarrow Q}{\Delta\rightarrow P\supset Q}\quad\quad\quad\f{\Delta\rightarrow P\quad\quad\quad\Delta\rightarrow P\supset Q}{\Delta\rightarrow Q}
\]
and as the last rule he has the rule of \emph{reductio ad absurdum}, in its strong, nonconstructive, version, which he calls also ``the rule of indirect proof'':
\[
\f{\Delta,\sim P\rightarrow Q\quad\quad\quad\Delta,\sim P\rightarrow \;\sim Q}{\Delta\rightarrow P}
\]

Then G\" odel claims that all the tautologies in the language with $\sim$, $\supset$ and $\rightarrow$ may be proved in this system, which should mean that if $P$ is a tautology, then the sequent $\rightarrow P$, with the empty left-hand side, is provable. If we want to introduce other connectives, G\" odel's way of dealing with that would be to have rules corresponding to the definitions of these connectives in terms of $\sim$ and $\supset$. (He had an analogous approach in the Notre Dame notes with systems of the Hilbert type.)

Then G\" odel states that this system has over the system of the Hilbert type that he considered previously the advantage that: ``All the axioms are really very simple and evident.'' Next he notes that the ``pseudo-paradoxical'' $q\rightarrow p\supset q$ and $\sim p\rightarrow p\supset q$ are provable in the system, although nobody can have any objection about what was assumed for the system, i.e.\ everybody would admit these assumptions ``if we interpret both the $\supset$ and the $\rightarrow$ to mean \emph{if\ldots then}''.  Then he proves these two pseudo-paradoxical sequents. The proof of the first, which uses the identity axiom, thinning and implication introduction, is straightforward. For the second, G\" odel presents the proof in the following form:
\begin{tabbing}
\hspace{7em}\=$p\rightarrow p$, \hspace{.7em}axiom\\*
\>$\sim p,p,\sim q\rightarrow p$, \hspace{.7em}by thinning\\
\>$\sim p\rightarrow \;\sim p$, \hspace{.7em}axiom\\
\>$\sim p,p,\sim q\rightarrow \;\sim p$, \hspace{.7em}by thinning\\
\>$\sim p,p\rightarrow q$, \hspace{.7em}by \emph{reductio ad absurdum}\\
\>$\sim p\rightarrow p\supset q$, \hspace{.7em}by implication introduction.
\end{tabbing}
After this proof, G\" odel says that he is sorry he has ``no time left to go into more details about this Gentzen system'', and proceeds to conclude the first part of the course, devoted to propositional logic.

\section{Comments on the natural deduction\\ system}
Ja\' skowski's natural deduction system of \cite{J34} has rules completely analogous to those of G\" odel's system: the rules for introducing and eliminating intuitionistic implication, which before the Second World War was called \emph{positive} implication, and strong \emph{reductio ad absurdum}. G\" odel's identity axioms are being taken care of in Ja\' skowski's system by posing hypotheses, while thinning is embodied in implication introduction.

G\" odel's inspiration for his natural deduction system need not however have come directly from \cite{J34}. A system of the Hilbert type for classical propositional logic in Frege's negation-implication language, which is due to {\L}ukasiewicz from \cite{LT30}, is obtained by simplifying the axioms from \cite{F79} (see the footnote concerning Frege after Theorem~5 in Section~2), and is mentioned in \cite{HA28} (Section I.10), which G\" odel knew, as we said in  Section~2. This system has, with modus ponens as a rule, axioms of the following forms (in G\" odel's notation):
\begin{tabbing}
\hspace{7em}\=$P\supset(Q\supset P)$,\\
\>$(P\supset(Q\supset R))\supset((P\supset Q)\supset(P\supset R))$,\\
\>$(\sim P\supset\;\sim Q)\supset(Q\supset P)$.
\end{tabbing}
(This system is also in \cite{Ch56}, Chapter II, Sections 20 and 27.) It is an exercise that G\" odel might well have done to show that strong contraposition, i.e.\ the third axiom-schema here, may be replaced by the axiom-schema of strong \emph{reductio ad absurdum}:
\begin{tabbing}
\hspace{7em}\=$(\sim P\supset Q)\supset((\sim P\supset\;\sim Q)\supset P)$.
\end{tabbing}
(The system L$_3$ of \cite{M64}, Section 1.6, amounts to the system with this replacement made; to show that strong contraposition yields strong \emph{reductio ad absurdum} one may first derive $\sim Q\supset(Q\supset\;\sim(P\supset P)$, and then use this to derive $\sim P\supset\;\sim(P\supset P)$ from the hypotheses $\sim P\supset Q$ and $\sim P\supset\;\sim Q$.)

It was well established at the time of G\" odel's Notre Dame course that, with modus ponens as a rule, the first two axiom-schemata above, stemming from \cite{F79} (Section 14), give a complete system for intuitionistic implication, which alternatively may be formulated by the rules for introducing and eliminating implication in natural deduction (this is established in \cite{J34}, Section~3, and was no doubt known to Gentzen). One passes from this basic implicational system, which falls short of classical implication (Peirce's law $((P\supset Q)\supset P)\supset P$ is missing), by adding to it a strong, nonconstructive, principle involving negation. This is either strong contraposition, or may be strong \emph{reductio ad absurdum}, as was Ja\' skowski's and G\" odel's choice, to make it fit better in the frame of natural deduction.

G\" odel was looking for a \emph{simple} system of propositional logic. Simple at the time when he was writing meant ``short'' or ``as economical as possible''. This economy was however sometimes spurious. Formulating propositional logic with the Sheffer stroke is a gain in economy paid with a high price in perspicuity. The pursuit of the least number of axioms, which was a popular subject before the Second World War, when it is not tied to an important independence problem, may also lead to exuberances. One can claim however that G\" odel was not looking for a spurious simplicity, but a genuine one.

With a slight change, he could have achieved another very simple presentation of propositional logic, in some sense simpler. Instead of negation he could have taken the propositional constant, i.e.\ nullary connective, $\bot$ as primitive, with the negation $\sim P$ being defined as $P\supset\bot$ (as suggested by Gentzen in \cite{G35}, Section II.5.2). Then his rule of \emph{reductio ad absurdum} could be replaced by the following version of this rule:
\[
\f{\Delta,P\supset\bot\rightarrow\bot}{\Delta\rightarrow P}
\]
which is shorter to state. Whether this other version is in general simpler could perhaps be questioned, since the constant $\bot$ is perhaps more puzzling than the connective of negation.

Another, more involved, change open to G\" odel, leading to another simple formulation of propositional logic, would be to take as primitive besides negation the connective of conjunction, for which he uses the symbol $.$ as in $P\, .\: Q$, instead of implication $\supset$. Then the two rules for $\supset$ could be replaced by the even simpler, and perfectly intuitive, following three rules for introducing and eliminating conjunction:
\[
\f{\Delta\rightarrow P\quad\quad\quad\Delta\rightarrow Q}{\Delta\rightarrow P\, .\:Q}\quad\quad\quad\f{\Delta\rightarrow P\, .\:Q}{\Delta\rightarrow P}\quad\quad\quad\f{\Delta\rightarrow P\, .\:Q}{\Delta\rightarrow Q}
\]
and the following simple form of the structural rule of cut:
\[
\f{R\rightarrow P\quad\quad\quad\Delta,P\rightarrow Q}{\Delta,R\rightarrow Q}
\]
which is ``invisible'' in formulations of natural deduction such as Gentzen's systems \emph{NJ} and \emph{NK} of \cite{G35}.
One defines $P\supset Q$ in terms of $\sim$ and $.$ by $\sim(P\, .\:\sim Q)$.

We prove $\sim\sim P\rightarrow P$ as in G\" odel's system, in the following manner:
\begin{tabbing}
\hspace{7em}\=$\sim P\rightarrow\;\sim P$, \hspace{.7em}axiom\\*
\>$\sim\sim P,\sim P\rightarrow\;\sim P$, \hspace{.7em}by thinning\\
\>$\sim\sim P\rightarrow\;\sim\sim P$, \hspace{.7em}axiom\\
\>$\sim\sim P,\sim P\rightarrow\;\sim\sim P$, \hspace{.7em}by thinning\\
\>$\sim\sim P\rightarrow P$, \hspace{.7em}by \emph{reductio ad absurdum}.
\end{tabbing}
We derive the weak, constructive, version of \emph{reductio ad absurdum}:
\[
\f{\Delta,P\rightarrow Q\quad\quad\quad\Delta,P\rightarrow \;\sim Q}{\Delta\rightarrow\;\sim P}
\]
in the following manner:
\begin{tabbing}
\hspace{7em}\=$\Delta,P\rightarrow Q$, \hspace{.7em}premise\\*
\>$\Delta,\sim\sim P\rightarrow Q$, \hspace{.7em}by simple cut with $\sim\sim P\rightarrow P$\\
\>$\Delta,P\rightarrow \;\sim Q$, \hspace{.7em}premise\\
\>$\Delta,\sim\sim P\rightarrow \;\sim Q$, \hspace{.7em}by simple cut with $\sim\sim P\rightarrow P$\\
\>$\Delta\rightarrow \;\sim P$, \hspace{.7em}by \emph{reductio ad absurdum}.
\end{tabbing}

For the defined $\supset$, G\" odel's rule for introducing implication is derived as follows:
\begin{tabbing}
\hspace{7em}\=$\Delta, P\rightarrow Q$, \hspace{.7em}premise\\*
\>$P\, .\sim Q\rightarrow P\, .\sim Q$, \hspace{.7em}axiom\\
\>$P\, .\sim Q\rightarrow P$, \hspace{.7em}by conjunction elimination\\
\>$\Delta, P\, .\sim Q\rightarrow Q$, \hspace{.7em}by simple cut\\
\>$P\, .\sim Q\rightarrow\;\sim Q$, \hspace{.7em}by conjunction elimination\\
\>$\Delta,P\, .\sim Q\rightarrow\;\sim Q$, \hspace{.7em}by thinning\\
\>$\Delta\rightarrow\;\sim(P\, .\sim Q)$, \hspace{.7em}by weak \emph{reductio ad absurdum}.
\end{tabbing}

For the defined $\supset$, G\" odel's rule for eliminating implication is derived as follows:
\begin{tabbing}
\hspace{7em}\=$\Delta\rightarrow P$, \hspace{.7em}premise\\*
\>$\Delta,\sim Q\rightarrow P$, \hspace{.7em}by thinning\\
\>$\sim Q\rightarrow\;\sim Q$, \hspace{.7em}axiom\\
\>$\Delta,\sim Q\rightarrow\;\sim Q$, \hspace{.7em}by thinning\\
\>$\Delta,\sim Q\rightarrow P\, .\sim Q$, \hspace{.7em}by conjunction introduction\\
\>$\Delta\rightarrow\;\sim(P\, .\sim Q)$, \hspace{.7em}premise\\
\>$\Delta,\sim Q\rightarrow\;\sim(P\, .\sim Q)$, \hspace{.7em}by thinning\\
\>$\Delta\rightarrow Q$, \hspace{.7em}by \emph{reductio ad absurdum}.
\end{tabbing}

The burden of cancellation, i.e.\ discharging, of hypotheses is now taken care of entirely by strong \emph{reductio ad absurdum}. Besides the form it has in G\" odel's system, this rule can also now be given the form it has in the following somewhat shorter version, involving conjunction:
\[
\f{\Delta,\sim P\rightarrow Q\, .\sim Q}{\Delta\rightarrow P}
\]

The peculiarity of this last negation-conjunction system is that although it is based on a strong, nonconstructive, version of \emph{reductio ad absurdum}, all its provable sequents of the form $\rightarrow P$, with an empty left-hand side, are intuitionistically correct, i.e.\ for them $P$ is intuitionistically correct. (This is not the case for all sequents $\Delta\rightarrow P$, as witnessed by the provable sequent ${\sim\sim p\rightarrow p}$, the proof of which is given above.) This consequence of Glivenko's theorem is drawn by G\" odel at the beginning of his paper on the double-negation translation \cite{Go33}. In the negation-conjunction language, intuitionistic and classical propositional logic do not differ if we take into account only categorical provability, but they do differ if we take also into account hypothetical proofs, i.e.\ deductions.

The simple form of cut assumed for this other system in addition to thinning, which G\" odel did assume, brings us to the topic of the other structural rules not mentioned by G\" odel, which are derivable in his system. One obtains permutation:
\[
\f{\Delta,P,Q,\Gamma\rightarrow R}{\Delta,Q,P,\Gamma\rightarrow R}
\]
and contraction:
\[
\f{\Delta,P,P\rightarrow Q}{\Delta,P\rightarrow Q}
\]
with the help of thinning and the rules for implication. The derivation of contraction (which is shorter) is made as follows:
\begin{tabbing}
\hspace{7em}\=$\Delta,P,P\rightarrow Q$, \hspace{.7em}premise\\*
\>$\Delta,P\rightarrow P\supset Q$, \hspace{.7em}by implication introduction\\
\>$\Delta,P\rightarrow P$, \hspace{.7em}from the axiom $P\rightarrow P$ by thinning\\
\>$\Delta,P\rightarrow Q$, \hspace{.7em}by implication elimination.
\end{tabbing}
Note that for this derivation it is essential that the implication elimination rule be given in G\" odel's \emph{additive} version, and not in the following, more usual, \emph{multiplicative} version:
\[
\f{\Gamma\rightarrow P\quad\quad\quad\Delta\rightarrow P\supset Q}{\Delta,\Gamma\rightarrow Q}
\]

The cut rule in the form:
\[
\f{\Gamma\rightarrow P\quad\quad\quad\Delta,P\rightarrow Q}{\Delta,\Gamma\rightarrow Q}
\]
is also derivable with the help of thinning and the rules for implication.

The rule:
\[
\f{\Delta\rightarrow Q\quad\quad\quad\Delta\rightarrow \;\sim Q}{\Delta\rightarrow P}
\]
has the following derivation (which is related to the proof of $\sim p,p\rightarrow q$ at the end of Section~4):
\begin{tabbing}
\hspace{7em}\=$\Delta\rightarrow Q$, \hspace{.7em}premise\\*
\>$\Delta,\sim P\rightarrow Q$, \hspace{.7em}by thinning\\
\>$\Delta\rightarrow \;\sim Q$, \hspace{.7em}premise\\
\>$\Delta, \sim P\rightarrow \;\sim Q$, \hspace{.7em}by thinning\\
\>$\Delta\rightarrow P$, \hspace{.7em}by \emph{reductio ad absurdum}.
\end{tabbing}

We have shown above how to derive in G\" odel's system the weak, constructive, version of \emph{reductio ad absurdum} from its strong version with the help of $\sim\sim P\rightarrow P$ (proved above) and cut. Once we have all that, it is not difficult to ascertain that G\" odel's system is indeed complete by comparing it with Gentzen's systems of \cite{G35}.

\section{Is G\" odel's system a natural deduction\\ system?}
We have called G\" odel's system a natural deduction system without further ado, as Gentzen's system of \cite{G36} (Section II.5) is always considered to be such a system. Some people might however be taken aback by this usage, for the simple reason that G\" odel and Gentzen rely on sequents. It is fine to qualify Ja\' skowski's system of \cite{J34} as natural deduction, but they would think that it is spurious to call G\" odel's system that way. This matter is on the verge of turning into a scholastic terminological dispute. (A balanced precise assessment of it, without however dealing with the coding mentioned at the end of this section, may be found in \cite{I10}, Chapter~2, and \cite{I14}; see also \cite{D99}, Sections 0.3.4-5.)

Once we introduce for systems of the Hilbert type the notion of proof from hypotheses, without which it is very cumbersome to work with such systems, we come very close to natural deduction. The difference then boils down to whether the deduction theorem is proved or assumed as the rule for implication introduction.

Natural deduction systems may however be formulated in the absence of implication, and what characterizes them in that case is the fact that proofs from hypotheses are built by being prolonged at the bottom. If these proofs are given the form of trees, we extend them by introducing new roots, and if they are given the form of sequences we extend them at their end, downwards.

Sequents may be used for the purpose of making explicit on what uncancelled hypotheses a formula in a proof from hypotheses depends. The fact that this is noted very precisely with sequents, and is not taken care of by a cancellation device, which is often more cumbersome, does not abolish the natural deduction character of proofs from hypotheses that we mentioned in the preceding paragraph. The elimination of logical constants with sequents tied to natural deduction proofs happens on the right-hand side of sequents, as well as their introduction. This is how the working at the side of roots, i.e.\ modifying proofs only at their ends by extending them downwards, is manifested. With sequents in Gentzen's systems of the \emph{LJ} and \emph{LK} kind of \cite{G35}, we modify the proofs also at the side of hypotheses, by extending the corresponding proofs in \emph{NJ} and \emph{NK} upwards. This is achieved with Gentzen's rules for introducing logical constant at the left-hand side of sequents. We have not gone very far from natural deduction with that. Only the building of proofs is done in a different manner.

Anyway, G\" odel's system based on sequents qualifies as a natural deduction system by building proofs in the right manner. Logical constants are introduced and eliminated at the right-hand side of sequents, and proofs are extended downwards.

A sequent system of the \emph{LJ} and \emph{LK} kind also differs from a usual natural deduction system by having structural rules explicit. In natural deduction these rules are implicit, hidden under other rules, invisible, because taken for granted. Thinning is usually hidden under implication introduction, and cut, i.e.\ composition of proofs, is taken for granted. G\" odel's system has moved in that respect slightly away from the customary natural deduction format. Thinning is explicit rather than hidden by assuming also the following version of implication introduction:
\[
\f{\Delta\rightarrow Q}{\Delta\rightarrow P\supset Q}
\]
G\" odel was right when he found this rule less natural than thinning, and he preferred to keep for implication introduction just the more convincing rule he assumed.

It is also not very important in the usual natural deduction systems that hypotheses are arranged in sequences. Both Gentzen and G\" odel presumably wish here to be precise, and this should be the reason why, as Gentzen, G\" odel takes sequents to be words of a language, i.e.\ finite sequences of symbols. It is hard to see in this choice something else than this wish for precision, although with it one moves also in the direction where more involved matters of general proof theory arise, and Gentzen's and G\" odel's choice is proved right (see \cite{D16a}).

Although he had sequences of formulae on the left-hand side of sequents, G\" odel eschewed assuming permutation and contraction, and hid them under other rules in the natural deduction style. For his purposes, which involved also a wish for economy, as well as simplicity, this was presumably the right choice. Permutation and contraction should be more bizarre for the untrained logician than thinning. Even logicians who are not proof-theorists might be annoyed by them.

We don't think that G\" odel's natural deduction system differs in an essential way from Ja\' skowski's system. We presume that, with sequents, G\" odel wished  to be at the same time more precise and more natural. Ja\'skowski is extremely precise, but his precision goes in another, mathematically less important, direction, though it is perhaps important for computing.  Mathematically, the two systems are essentially the same.

We don't think that it makes a crucial difference
that Ja\' skowski arranges his proofs in boxes, explicit or implicitly marked with numbers. It is also not a crucial difference that Gentzen has his proofs in the systems \emph{NJ} and \emph{NK} of \cite{G35} built as trees, while G\" odel makes sequences as in systems of the Hilbert type, and notes dependencies in the margin, as at the end of Section~4, and later in our paper.

The tree form arises later in general proof theory when proofs are coded by terms---either typed lambda terms as in the Curry-Howard correspondence, or arrow terms of categorial proof theory (see \cite{D10} and \cite{D12}). This tree form is the form of terms where rules of inference appear as operations on terms. As it happens with Gentzen in other matters too, by pursuing his instinct of simplicity and precision he had a fine premonition of things that will become important with more advanced questions of general proof theory (see, for example, \cite{DP04}, Section 11.1).  However, without the coding of proofs (which becomes important if we are looking for criteria of identity of proofs; see \cite{DP04}, Chapter~1, and \cite{D10}), the form natural deduction proofs may take is more varied.

It is usually taken that natural deduction is characterized by assuming mainly rules and few, if any, axioms. In this sense, G\" odel's system is certainly a natural deduction system. Its only axiom ${P\rightarrow P}$ corresponds to assuming hypotheses in the format of Ja\' skowski or Gentzen's NJ and NK.

This matter of putting the burden on rules is however something natural deduction shares with sequent systems in general. It is tied to a shift from the categorical to the hypothetical in logical investigations. With natural deduction and sequents it is not only theorems that are important, but also deductions, and deductions with all  hypotheses cancelled are not privileged. If deductions are conceived only as consequence relations, we are still not very far from theorems. They should be conceived as arrows in a graph such that there may be more than one deduction between a premise and a conclusion in order to make us enter into another realm (see \cite{D10}, \cite{D12} and \cite{D16}).

After our considerations in the last part of Section~2 of \cite{DA16}, it is difficult to say that G\" odel had a premonition of that realm, as it seems Gentzen had. His presentation of natural deduction with the help of sequents was however a step in the right direction. With sequents one takes this road towards the investigation of deduction in a clearer and more decisive manner than in natural deduction without them.

\end{document}